\documentclass[11pt]{amsart}
\usepackage{amssymb}
\usepackage{amsmath}
\usepackage{amsthm}
\usepackage{xcolor}
\usepackage{pstricks,pst-node}

\newtheorem{theorem}{Theorem}[section]
\newtheorem{lemma}[theorem]{Lemma}
\newtheorem{definition}[theorem]{Definition}
\newtheorem{example}[theorem]{Example}
\newtheorem{proposition}[theorem]{Proposition}
\newtheorem{corollary}[theorem]{Corollary}
\newtheorem{remark}[theorem]{Remark}

\def\F{\mathcal{F}}
\def\p{\mathcal{P}}
\def\v{\overset\vee}

\def\dist{\mathrm{dist}}

\begin{document}
\title{Envelopes of holomorphy and extension of functions of bounded type.}
\author{Daniel Carando}
\thanks{This work was partially supported by CONICET, by ANPCyT PICT 05 17-33042, ANPCyT PICT 06 00587 and UBACyT Grant X038. The second author was also supported by a doctoral fellowship from CONICET}
\author{Santiago Muro}
\address{Departamento de Matem\'{a}tica - Pab I,
Facultad de Cs. Exactas y Naturales, Universidad de Buenos Aires,
(1428) Buenos Aires, Argentina, and CONICET}
\email{dcarando@dm.uba.ar} \email{smuro@dm.uba.ar}

\date{}

\begin{abstract}
We study the extension of holomorphic functions of bounded type defined on an open subset of a Banach space, to larger domains. For this, we first characterize the envelope of holomorphy of a Riemann domain over a Banach space, with respect to the algebra of bounded type holomorphic functions, in terms of the spectrum of the algebra. We then give a simple description of the envelopes of balanced open sets and relate the concepts of domain of holomorphy and polynomial convexity. We show that for bounded balanced sets, extensions to the envelope are always of bounded type, and that this does not necessarily hold for unbounded sets, answering a question posed by Hirschowitz in 1972.  We also consider extensions to open subsets of the bidual, present some Banach-Stone type results and show some properties of the spectrum when the domain is the unit ball of $\ell_p$.
\end{abstract}


\maketitle

\section*{Introduction}
This work was motivated by the following question: given an open subset $U$ of a complex Banach space, which is the largest open set containing $U$ to which holomorphic functions of bounded type on $U$ extend uniquely? As it could be expected, to properly pose and study the problem, we must expand our investigations to the Riemann domain framework. Our problem translates, then, to the characterization of the envelope of holomorphy of a Riemann domain modeled on a Banach space, with respect to the algebra of analytic funcions of bounded type. Loosely speaking, if $X$ is a Riemann domain over the Banach space $E$, the $H_b$-envelope of holomorphy of $X$ is the largest Riemann domain (over $E$) ``containing $X$'' to which every holomorphic function of bounded type on $X$ has a unique holomorphic extension.
We show in Theorem~\ref{envelopes coincide} that the $H_b$-envelope of holomorphy of $X$ can be identified with a subset of the spectrum of $H_b(X)$, the algebra of all holomorphic functions on $X$ of bounded type.

When we turn back to our original motivation and start with an open set $U\subset E$, we want to find out conditions on $U$ that ensures its envelope to be also an open subset of $E$ (in this case, the envelope is said to be schlicht or univalent). In this context, we will show in Theorem~\ref{envoltura H_b balanceados} that, if $U$ is a  balanced open subset, then its envelope of holomorphy is univalent and has a simple description in terms of the polynomially convex hull of $U$.
If $U$ is also bounded, extensions to the envelope are of bounded type (Theorem~\ref{extension acotados}). However, we show in Example~\ref{ejemplo}  an unbounded balanced open subset of $c_0$ for which extensions to the envelope are not necessarily of bounded type, answering a question of Hirschowitz~\cite[Remarque 1.8]{Hir72}.
We also relate the polynomial convexity with the property of being a $H_b$-domain of holomorphy, showing that these two concepts coincide for bounded balanced open sets (Theorem~\ref{extension acotados}).

In Section~\ref{seccion-bidual} we consider extensions from $U\subset E$ to open subsets of the bidual of $E$,
with particular interest on the Aron-Berner extension~\cite{AroBer78}. For a balanced subset $U$ of a symmetrically regular Banach space $E$, we describe in Proposition~\ref{extension E segunda mas} the largest open subset of $E''$ to which there exist Aron-Berner extensions of functions in $H_b(U)$. This set could be seen as the envelope of $U$ in the sense of Dineen and Venkova's work~\cite{DinVen04}.

In section 4 we consider Banach spaces for which finite type polynomials are dense in $H_b(E)$. When they are also reflexive, they are called Tsirelson-like spaces following~\cite{Vie07}.  We  show in Theorem~\ref{tsirelson_reciproco} that Tsirelson-like spaces are precisely the spaces where the holomorphic convexity of a balanced open set $U$ is equivalent to all the elements of the spectrum being evaluations at points in $U$, extending some results of \cite{Muj01} and \cite{Vie07}. This means that Tsirelson-like spaces are the only spaces that behave as in the several complex variables theory. We characterize the density of finite type polynomials in terms of the image of the spectrum of $H_b(U)$ by its canonical projection on $E''$. We also give a Banach-Stone type result (Theorem~\ref{Banach-Stone}) which improves some results in \cite{Vie07} and \cite{CarGarMae05}.

In the last section we present some properties of $M_b(U)$, the spectrum of $H_b(U)$, somehow extending the study of \cite{AroGalGarMae96} and \cite{CarGarMae05}. In the case $U=E$, it was shown in \cite[Section 6.3]{Din99} that bounded type entire functions extend to holomorphic functions on the spectrum of $H_b(E)$ which are of bounded type on each connected component. We prove, in contrast, that in most cases there are polynomials whose extensions are not of bounded type on the whole Riemann domain $M_b(E)$. Then we concentrate in the case $U=B_{\ell_p}$ to show that the structure of the spectrum in not what one may expect from the case $U=E$, with $E$ a symmetrically regular Banach space. In the latter case, $M_b(E)$ is the disjoint union of copies of $E''$. However, we show that $M_b(B_{\ell_p})$ is not a disjoint union of ``unit balls".
For $p\in\mathbb N$, we also define a distinguished spectrum to where the canonical extensions are of bounded type and which turns out to be a $H_b$-domain of holomorphy.

\medskip

In the theory of several complex variables, where every holomorphic function is of bounded type, the envelope of holomorphy
can be described in terms of
the spectrum of the algebra $H(X)$ of all analytic functions on $X$.
This idea stems from Bishop, who introduced an analytic structure in the spectrum that makes it a Riemann domain \cite{Bis63}. In the infinite dimensional setting, most of the study was done for the space of all holomorphic functions on open subsets of Banach and more general locally convex spaces (see, for example, \cite{Mat72,Muj86,Muj87,Sch83}). The study of the spectrum of the algebra of holomorphic functions of bounded type on a Banach space was initiated in 1991 by Aron, Cole and Gamelin in their seminal article \cite{AroColGam91}.
Recall that a Banach space $E$ is said to be (symmetrically) regular if every continuous (symmetric) linear mapping $T:E \to E'$ is weakly
compact (an operator $T:E\to E'$ is symmetric if $Tx_1(x_2) = Tx_2(x_1)$ for all $x_1,x_2 \in E$). In \cite[Corollary 2.2]{AroGalGarMae96} Aron, Galindo, Garc\'{\i}a and
Maestre gave ${M}_{b}(U)$ a structure of Riemann analytic
manifold modeled on $E''$, for $U$ an open subset of a symmetrically regular space $E$.
For the case $U=E$, ${M}_{b}(E)$ can be viewed as the
disjoint union of analytic copies of $E''$, these copies being
the connected components of ${M}_{b}(E)$. In
\cite[Section 6.3]{Din99}, there is an elegant exposition of many of
these results. The study of the spectrum of the algebra of the space
of holomorphic functions of bounded type was continued in
\cite{CarGarMae05}. The analytic structure of $M_b(X)$ for $X$ a Riemann domain over a symmetrically regular Banach space $E$ was presented by Dineen and Venkova in \cite{DinVen04}. 

\bigskip

Throughout this paper $E$ and $F$ will be complex Banach
spaces. We denote by $\p(^nE)$ the Banach space of all
continuous $n$-homogeneous polynomials from $E$ to $\mathbb C$, and by $\p(E)$ the class of all continuous polynomials.
If $U$ is an open subset of $E$, $H(U)$ denotes the space of all holomorphic functions on $U$. A holomorphic function $f\colon U\to \mathbb C$ is \emph{of bounded type} if it is bounded on \emph{$U$-bounded sets} (i.e., bounded subsets that are bounded away from the boundary of $U$). We denote by $H_b(U)$ the space of all analytic functions of bounded type, which is a Fr\'{e}chet space when endowed with the topology of uniform convergence on $U$-bounded sets. It is known that for balanced open sets $U$, polynomials are dense in $H_b(U)$ (see for example \cite[Theorem 7.11]{Muj86}).

We refer the reader to the already mentioned articles \cite{AroColGam91,AroGalGarMae96,DinVen04} for a description of $M_b(U)$ and $M_b(X)$ and their  analytic structure (see also~\cite{CarGarMae05}), and to the books by Dineen \cite{Din99} and Mujica \cite{Muj86} for a more extensive treatment of infinite dimensional holomorphy.

\section{The $H_b$-envelope of holomorphy}

A Riemann domain  $(X,p)$ over the Banach space $E$ is a Hausdorff topological $X$ space with a local homomorphism $p:X\to E$. For each $x\in X$ we define $\dist_X(x)$ as the supremum of all  $r>0$ for which there exists a neighborhood of $x$ homeomorphic via $p$ to $B(p(x),r)$.
We say that  $A\subset X$ is $X$-bounded if $p(A)$ is bounded and \begin{equation}\label{defi-distancia}\dist_X(A):=\inf\{\dist_X(x):\, x\in A\}\end{equation} is positive.  The definition of $H_b(X)$ is now clear: bounded type holomorphic functions on $X$ are those which are bounded on $X$-bounded sets.
The space $H_b(X)$ is a Fr\'echet algebra when it is considered with the topology of uniform convergence on $X$-bounded sets. By a \textit{fundamental sequence of $X$-bounded sets} we  mean a sequence $\{A_n\}_n$ of $X$-bounded subsets such that if $B$ is another $X$-bounded subset, then there exists $n_0$ with $B\subset A_{n_0}$. A typical fundamental sequence of $X$-bounded sets is given by $(X_n)_n$ with \begin{equation}\label{typical fundamental sequence}X_n:=\{x\in X:\, \dist_X(x)\ge\frac1{n} \textrm{ and }\|p(x)\|\le n\}.\end{equation}

We will denote by $M_b(X)$ the spectrum of the algebra $H_b(X)$, that is, the set of all non-zero continuous, linear and multiplicative functionals on $H_b(X)$. Thus, for each $\varphi \in M_b(X)$ there exists an $X$-bounded set $B$ such that $|\varphi(f)|\le \sup_{x\in B}|f(x)|$, for all $f\in H_b(X)$. In this case, we will write $\varphi\prec B$. We also set $$\|f\|_B:=\sup_{x\in B}|f(x)|.$$

For $X=U\subset E$ an open subset, we define an application $\pi:M_b(U)\to E''$ by $\pi(\varphi)=\varphi_{|_{E'}}$. If $E$ is symmetrically regular, this mapping $\pi$ provides the local homeomorphism that makes $M_b(U)$ a Riemann domain over $E''$ \cite{AroGalGarMae96}.
For a general Riemann domain $(X,p)$ modeled over a symmetrically regular Banach space $E$, the mapping $\pi:M_b(X)\to E''$ is defined by $\pi(\varphi)(\gamma)=\varphi(\gamma\circ p)$, and the analytic structure is analogous \cite{DinVen04}. Functions in $H_b(X)$ naturally extend to $M_b(X)$ by the Gelfand transform, and it is shown in \cite{AroGalGarMae96} and \cite{DinVen04} that this extension is analytic.
Symmetric regularity is necessary for the analytic structure of  $M_b(X)$ to work \cite[Proposition 2.3]{AroGalGarMae96}, because one has to deal with Aron-Berner extensions. However, if we restrict ourselves to  $\pi^{-1}(E)\subset M_b(X)$, we can do without symmetric regularity. To see this,  we first recall that if $f$ is a holomorphic function on $X$, its differential at $y\in X$ is given by
$$d^nf(y)= d^n[ f\circ(p_{|_{V_{y}}})^{-1}]\big(p(y)\big),$$ where $V_y$ is some neighbourhood of $y$ on which $p$ is homeomorphic. Now, fix
$\varphi\in\pi^{-1}(E)$ and $\delta<1/\dist_X(\varphi)$. We can define for each $x\in E$ with  $\|x\|<\delta$
\begin{equation}\label{fi supra x}
{\varphi^x(f)=\sum_{n=0}^\infty\varphi\Big(\frac{d^nf(\cdot)}{n!}(x)\Big).}
\end{equation}
Just as in \cite{AroGalGarMae96} or \cite{DinVen04}, the sets $\{\varphi^x:\; \|x\|<\delta\}$ (for $\varphi\in\pi^{-1}(E)$ and $\delta<\dist_X(\varphi)$) form a basis of a Hausdorff topology for $\pi^{-1}(E)$, and \begin{equation}\label{pi restringida}\pi|_{\pi^{-1}(E)}:\pi^{-1}(E)\to E \end{equation} is a local homeomorphism on each of these sets. This endows $\pi^{-1}(E)$ with an analytic structure over $E$, for which the Gelfand transform of any function in $H_b(X)$ is analytic. So we have.
\begin{lemma}\label{pi a la menos 1 dom de riemann}
Let $(X,p)$ be a Riemann domain spread over a Banach space $E$  and let $\pi:M_b(X)\to E''$ be defined as above. Then $(\pi^{-1}(E),\pi)$ is a Riemann domain spread over $E$ (here we consider the restriction of $\pi$ as in Equation~\ref{pi restringida} ). Also, any  $f\in H_b(X)$ extends to an analytic function on $\pi^{-1}(E)$ via the (restriction of) the Gelfand transform.
\end{lemma}

\bigskip

We recall the definition of extension morphism and envelope of holomorphy for a family of holomorphic functions (see, for example, \cite[Chapter XIII]{Muj86}).
Let $(X,p)$ and  $(Y,q)$ be Riemann domains spread over a Banach space $E$. A \emph{morphism} from $X$ to $Y$ is a continuous mapping $\tau:X\to Y$ such that $q\circ \tau=p$. If $\mathcal F$ is a subset of $H(X)$, then a morphism $\tau:X\to Y$ is said to be an \emph{$\mathcal F$-extension} of $X$ if for each $f\in \mathcal F$ there is a unique $\tilde f\in H(Y)$ satisfying $\tilde f\circ \tau=f$, so that the following diagram commutes:
$$
\begin{psmatrix}[colsep=1cm,rowsep=1.2cm]
  X        &    &  Y \\
          &   E  &
\psset{arrows=->,labelsep=3pt,nodesep=3pt}
\ncline{1,1}{1,3}^{\tau}
\ncline{1,1}{2,2}<{p}
\ncline{1,3}{2,2}>{q}
\end{psmatrix}
$$
A morphism $\tau:X\to Y$ is said to be an \emph{$\mathcal F$-envelope of holomorphy} of $X$ if $\tau$ is an $\mathcal F$-extension of $X$ and if for each $\mathcal F$-extension $\nu:X\to Z$ of $X$, there is a morphism $\mu:Z\to Y$ such that $\mu\circ\nu=\tau$:
$$
\begin{psmatrix}[colsep=1.8cm,rowsep=1.2cm]
          &                 &  Y \\
  X       &   Z             & \\
          &   E             &
\psset{arrows=->,labelsep=3pt,nodesep=3pt}
\ncarc[arcangle=20]{2,1}{1,3}^{\tau}
\ncline[linestyle=dashed,dash=3pt 2pt]{2,2}{1,3}<{\mu} 
\ncarc[arcangle=20]{1,3}{3,2}>{q}
\ncline{2,1}{2,2}^{\nu}
\ncline{2,2}{3,2}<{}
\ncline{2,1}{3,2}<{p}
\end{psmatrix}
$$
Regarding holomorphic functions of bounded type, the $H_b$-envelope of holomorphy was constructed  by Hirschowitz in \cite{Hir72} by means of germs. For general families of functions $\mathcal F$, the existence of the $\mathcal F$-envelope of holomorphy can be seen in \cite[Chapter XIII]{Muj86}. We say that a Riemann domain $(X,p)$ is a \emph{$H_b$-domain of holomorphy} if the identity on $X$ is an $H_b$-envelope of  holomorphy. Loosely speaking, this is to say that $X$ coincides with its $H_b$-envelope of holomorphy.

Our concept of $H_b$-extension morphism is different from that
introduced by Dineen and Venkova in \cite{DinVen04}. The main difference is that in our case, the envelope of a Riemann domain over $E$ is also modeled on $E$, while theirs is modeled on $E''$ (just as the spectrum \cite{AroGalGarMae96}).

Now we are ready to give the characterization of the $H_b$-envelope of holomorphy, which is very similar to that of several complex variables, especially if  $E$ is reflexive.

\begin{theorem}\label{envelopes coincide}
Let $(X,p)$ be a connected Riemann domain spread over a Banach space $E$ and let $Z$ be the connected component of $\pi^{-1}(E)\subset M_b(X)$ which intersects $\delta(X)$. Then $\delta:(X,p)\to(Z,\pi)$, $\delta(x)=\delta_x$ is the $H_b$-envelope of $X$.
\end{theorem}
\begin{proof}
Denote by $\tau:(X,p)\to (Y,q)$ the $H_b$-envelope of $X$.
As $\delta:(X,p)\to(Z,\pi)$  is an $H_b$-extension, there exists a morphism $\mu:Z \to  Y$ such that $q=\mu\circ\delta$.
For each $f\in H_b(X)$ we have a unique $\tilde f\in H(Y)$ with $\tilde f\circ \tau=f$. Since $(Z,\pi)$ is an $H_b$-extension of $(X,p)$, the function $\tilde f\circ \mu\in H(Z)$ must be the restriction to $Z$ of the canonical extension of $f$ to $M_b(X)$. Therefore,  $\varphi (f)= \tilde f(\mu(\varphi))$ for all $f$ and then $\mu(\varphi)$ uniquely determines $\varphi$. This means that $\mu$ is injective.

Now we see that $\mu$ is onto.
Since $\mu$ is a morphism, $\mu(Z)$ is open in $Y$.
Suppose that there exists $y\in\overline{\mu(Z)}\setminus\mu(Z)$. For a fundamental sequence of $X$-bounded sets  $(X_n)_n$  as in~(\ref{typical fundamental sequence}), we define $W_n=\{\varphi\in Z:\, \varphi\prec X_n\}$.
From~\cite{AroGalGarMae96} (see also \cite[Proposition 1.5]{DinVen04}) we have $\dist_X(W_n)\ge1/n$ and therefore we can get  a sequence of natural numbers $(n_k)_k\in \mathbb N$ and, for each $k$, a homomorphism  $\varphi_k\in W_{n_{k+1}}\setminus W_{n_k}$ so that the sequence $\mu(\varphi_k)$ converges to $y$ in $Y$. We can now choose functions $f_k\in H_b(X)$ for which $\|f_k\|_{X_{n_k}}<1/2^k$ and $$|\varphi_k(f_k)|>k+\sum_{j=1}^{k-1}| \varphi_k(f_j)|.$$ The series $\sum_{j=1}^{\infty}f_j$ converges to some $f\in H_b(X)$ and, since $\varphi_k$ is a continuous homomorphism,  we have $$|\varphi_k(f)| = \Big|\sum_{j=1}^{\infty} \varphi_k(f_j)\Big| \ge| \varphi_k(f_k)|-\Big|\sum_{j=1}^{k-1} \varphi_k(f_j)\Big|-\Big|\sum_{j=k+1}^{\infty} \varphi_k(f_j) \Big|>k-1.
$$
This means that $|\varphi_k(f)|\to\infty$ as $k\to \infty$. But if we take $\tilde f\in H(Y)$ satisfying $f=\tilde f\circ \tau$,  we also have $\varphi_k(f)=\tilde f(\mu(\varphi_k))$, which tends to $\tilde f(y)$, a contradiction. Thus, $\mu(Z)$ is closed in $Y$.
Since $Y$ is connected, we conclude that $\mu(Z)=Y$.
$\square$ \end{proof}

Theorem \ref{envelopes coincide} states in particular that the $H_b$-envelope is part of the spectrum. In other words, evaluations at elements of the $H_b$-envelope are always continuous. This also happens for evaluations at elements on any other $H_b$-extension. Indeed,  we can proceed as in the beginning of the previous proof to show that the evaluation on an element of any extension coincides with the evaluation on some element of the $H_b$-envelope, which is continuous.

\section{Envelopes of open subsets of a Banach space}

In this section we restrict ourselves to open subsets of a Banach space $E$. In order to give a more precise and concrete description of the $H_b$-envelope of an open set $U\subset E$, we first define certain open sets which contain $U$ and to which, under some conditions, functions in $H_b(U)$ extend.

Let $U\subset E$ be an open set and $\F$ be a set of functions defined on $E$ (mainly,  $\F$ will be $H_b(E)$, $\p(E)$ or $\p(^nE)$ for some $n\in\mathbb N$). For
 $A$ a $U$-bounded set, we define the \emph{$\F$-hull} of $A$ as the set 
$$
\widehat A_\F=\{x\in E:\; |f(x)|\le \|f\|_A \;\; \textrm{for every } f\in\F\}.
$$
For $\F=H_b(U)$ (or any set of functions defined on $U$), the definition of $\widehat A_\F$ is analogous, just taking $x\in U$
instead of $x\in E$ in the set above.
Given a fundamental sequence of $U$-bounded sets $(U_n)_n$, such as the one constructed in~(\ref{typical fundamental sequence}), we define the \emph{$\F$-hull} of $U$ by
$$
\widehat U_\F:=\bigcup_{n\in\mathbb N}\Big(\widehat U_n\Big)_\F.
$$

An open set $U$ is \emph{$\F$-convex} if $\widehat A_{\F}$ is $U$-bounded for every $U$-bounded set $A\subset U$.
Our definition of $\F$-convex set coincides with the notion of strongly $\F$-convex set investigated by Vieira in \cite{Vie07}.
She proved that  $\widehat A_{\p(E)}=\widehat A_{H_b(E)}$ for each bounded set $A$, and as a consequence  $U$ is $\p(E)$-convex if and only if it is $H_b(E)$-convex . Also, it is easy to see that if $U$ is $\p(E)$-convex, then it is also $H_b(U)$-convex.
Whenever $U$ is balanced,  the $\p(E)$-convexity of $U$  is equivalent to its $H_b(U)$-convexity \cite[Proposition 1.5]{Vie07}.

We say that a point $x\in E$ is \emph{an evaluation for $H_b(U)$} if there is some $\varphi\in M_b(U)$ such that $f(x)=\varphi(f)$ for every $f\in H_b(E)$. If $H_b(E)$ is dense in $H_b(U)$, then $\varphi$ is uniquely determined by $x$, and we denote it by  $\delta_x$. The set of all evaluation points for $H_b(U)$ will be denoted by $\v U$. So we have the following.

\begin{lemma}\label{propiedades de los techitos}
Let $U$ be an open subset of the Banach space $E$. For $\F=\p(E)$ or $\F=H_b(E)$ we have:
\begin{enumerate}
\item[\rm(a)] both $\v U$ and $\widehat U_\F$ are open subsets of $E$, and  $U\subset \v U\subset \widehat U_\F$;
\item[\rm(b)] the set $U$ is $\F$-convex if and only if $U= \v U= \widehat U_\F$;
\item[\rm(c)] if $U$ is balanced, then $\v U= \widehat U_\F$, and they identify with the connected component of $\pi^{-1}(E)$ which intersects $U$ (where $\pi:M_b(U)\to E''$ is the local homomorphism).
\end{enumerate}
\end{lemma}
\begin{proof}
(a) It is known that $\widehat U_\F$ is open. To see that so is $\v U$, we take $x\in \v U$ and choose $\varphi\in M_b(U)$ such that $\varphi(f)=f(x)$ for every entire function of bounded type $f$. In particular, $\pi(\varphi)=x$ and then $\varphi$ actually belongs to $\pi^{-1}(E)$.
We can take $\delta>0$ such that the homomorphisms  $\varphi^y$ defined as in (\ref{fi supra x}) belong to $\pi^{-1}(E)$ for every $y\in B_E(0,\delta)$. Since for every $f\in H_b(E)$ we have  $\varphi^y(f)=f(x+y)$, we conclude that $x+y$ is in $\v U$ for any $y\in B_E(0,\delta)$.

In order to show the inclusion $\v U\subset \widehat U_\F$, for $z\notin \widehat U_\F$ we take functions $f_n\in\F$ such that $f_n(z)=1$ and $\|f_n\|_{U_n}\to 0$.
Thus, the evaluation at $z$ of bounded type entire functions is not continuous with the topology induced by $H_b(U)$, and therefore $z$ does not belong to $\v U$.

(b) The ``only if'' part is a consequence of the definitions and of (a). The ``if'' part  follows from \cite[Lemma 1.3]{Vie07}.

(c)
If $x\in\widehat U_\F$, then there exists $n\in\mathbb N$ such that $|f(x)|\le\|f\|_{U_n}$ for every $f\in \p(E)$. The homomorphism $\delta_x$ is then bounded on $\p(E)$ with the toplogy induced by $H_b(U)$. Since $U$ is balanced, $\p(E)$ is dense in $H_b(U)$, so $\delta_x$  extends to an element of $M_b(U)$. This shows that $\widehat U_\F$ coincides with $\v U$ (the reverse inclusion was given in (a)) and that it identifies with a subset of $\pi^{-1}(E) \subset M_b(U)$.
Moreover, by \cite[Lemma 1.4]{Vie07}, $\widehat U_{\p(E)}$ is balanced and hence connected. Let us see that it is a connected component of $\pi^{-1}(E)$.
We  denote by $\pi_E:M_b(E)\to E''$ the natural projection associated to the spectrum of $H_b(E)$, to distinguish it from $\pi:M_b(U)\to E''$. By Lemma~\ref{pi a la menos 1 dom de riemann}, $(\pi_E^{-1}(E), \pi_E)$
is a Riemann domain, and it can be seen as a disjoint union of copies of $E$, just as in \cite[Corollary 2.5]{AroGalGarMae96} and the comments previous to it. Since $H_b(E)$ is dense in $H_b(U)$, the subset $\pi^{-1}(E)$ of $M_b(U)$ (which is also a Riemann domain over $E$ by Lemma~\ref{pi a la menos 1 dom de riemann})
may be embedded in $(\pi_E^{-1}(E), \pi_E)$. Thus, the
connected component of $\pi^{-1}(E)$ which contains $U$ must be a
subset of $E$.
Now, if we take  $z\in E\setminus \v U$, there cannot be a $\phi\in M_b(U)$ which coincides with $\delta_z$ on entire functions (or on polynomials). So the connected component of $\pi^{-1}(E)$ which intersects $U$ is $\v U$.
$\square$ \end{proof}

We are ready to state our characterization of the $H_b$-envelope of holomorphy of an open balanced set.

\begin{theorem}\label{envoltura H_b balanceados}
Let $U$ be an open balanced subset of a Banach space $E$.
Then $\widehat U_{\p(E)}$ is the $H_b$-envelope of $U$. Moreover, any  $f\in H_b(U)$ extends to a holomorphic function $\tilde f$ on $\widehat U_{\p(E)}$ which is bounded on $\widehat A_{\p(E)} $ for every $U$-bounded set $A$. Also, for $z\in \widehat A_{\p(E)}$ we have $\delta_z\prec A$.
\end{theorem}
\begin{proof}
A direct combination of Theorem~\ref{envelopes coincide} and Lemma~\ref{propiedades de los techitos}(c) gives the first assertion, and the fact that any $f\in H_b(U)$ extends to a holomorphic function $\tilde f$ on $\widehat U_{\p(E)}$.
The remaining assertions are a consequence of the denseness of $\p(E)$ in $H_b(U)$, and are essentially proved in Lemma~\ref{propiedades de los techitos}.
$\square$ \end{proof}

It should be pointed out that we cannot expect to extend all the functions of $H_b(U)$ to connected sets with points outside $\v U$. Indeed, suppose $V\supset U$ is any  connected open set such that the inclusion $U\hookrightarrow V$ is an $H_b$-extension. If $z\in V$ then $\delta_z$ belongs to $M_b(U)$ by Theorem~\ref{envelopes coincide} and the comments after it. Moreover, we have $\delta_z(f)=f(z)$ for  every entire function $f$, so we conclude that $z$ must belong to $\v U$. By Lemma~\ref{propiedades de los techitos}(a),  we cannot either extend every function in $H_b(U)$ outside $\widehat U_{\p(E)}$.

It is natural to ask if the extension of a holomorphic function of bounded type on a balanced open set $U$ to $\widehat U_{\p(E)}$ is necessarily of bounded type. This is the case when $U$ is also bounded. To prove this we will use the following Lemma which states that the polynomial hull of a balanced set coincides with the intersection of its homogeneous polynomial hulls. This was noticed, for example, in \cite{Sic85} for balanced sets in $\mathbb C^n$.

\begin{lemma}\label{envoltura homogenea}
Let $V\subset E$ be an open balanced set. Then $$\widehat V_{\p(E)}=\bigcap_{n\in\mathbb N}\widehat V_{\p(^nE)}.$$
\end{lemma}
\begin{proof}
We only need to prove that $\bigcap_{n\in\mathbb N}\widehat V_{\p(^nE)}\subset\widehat V_{\p(E)}$, since the other inclusion is clearly true for every open set.
Suppose $z$ belongs to $\widehat V_{\p(^nE)}$ for all $n\in \mathbb N$ and let $P\in\mathcal P(E)$ be a polynomial of degree $k$. For any $n\in\mathbb N$ we can write $P^n=Q_0+\dots+Q_{nk}$, with $Q_j\in\mathcal P(^jE)$. By  Cauchy's inequalities, $\|Q_j\|_V\le\|P^n\|_V$, and thus $$|P^n(z)|=|\sum_{j=0}^{nk}Q_j(z)|\le\sum_{j=0}^{nk}\|Q_j\|_V\le\sum_{j=0}^{nk}\|P^n\|_V=(nk+1)\|P\|_V^n.$$ We then have  $|P(z)|\le(nk+1)^{\frac1{n}}\|P\|_V$ for every $n\in\mathbb N$, which implies that $|P(z)|\le\|P\|_V$.~$\square$ \end{proof}

\begin{theorem}\label{extension acotados}
If $U\subset E$ is a bounded and balanced open set, then every function in $H_b(U)$ can be extended to a holomorphic function of bounded type in $\widehat U_{\p(E)}$.
Moreover, $\widehat U_{\p(E)}$ is a $H_b$-domain of holomorphy. Hence, $U$ is a  $H_b$-domain of holomorphy if and only if $U=\widehat U_{\p(E)}$.
\end{theorem}
\begin{proof} Since $U$ is a bounded balanced set, for any sequence of
positive numbers $r_n \uparrow 1$, the sets $(r_n U)_{n\in\mathbb N}$
form a fundamental system of $U$-bounded sets. But $\widehat U_{\p(E)}$ is
also bounded and balanced, so $(r_n \widehat U_{\p(E)})_{n\in\mathbb N}$ is
a fundamental system of $\widehat U_{\p(E)}$-bounded sets.
Therefore, to see that extensions are of bounded type on $\widehat
U_{\p(E)}$,
it is enough to show that $r
\widehat U_{\p(E)}$ is contained in $(r U)^{^\wedge}_{\p(E)}$ for each $r<1$ and
use
Theorem~\ref{envoltura H_b balanceados}. Now, for $x\in r\widehat
U_{\p(E)}$ and any $Q\in \mathcal P(^jE)$ ($j\in \mathbb N$)
we have
$$|Q(x)|= r^j \big| Q\big(x/r)\big| \le r^j \|Q\|_U= \|Q\|_{r U}.$$
This means that $x$ belongs to $(rU)^{^\wedge}_{\mathcal P_j(E)}$ for any
$j\in \mathbb N$, so  by Lemma \ref{envoltura homogenea} we can conclude that
$x\in(rU)^{^\wedge}_{\p(E)}$.

To prove the second statement, take $z\in \widehat U_{\p(E)}$ and let $A$
be a $U$-bounded set such that $z\in \widehat A_{\p(E)}$. By Theorem~\ref{envoltura H_b balanceados} we have  $\delta_z\prec A$.
If we take $r<\dist(A,E\setminus U)$, then  the homomorphism
$(\delta_z)^y$ given as in (\ref{fi supra x}) is well defined and continuous on $H_b(U)$ for each $y\in rB_E$. This means that $B_{E}(z,r)$ is contained in
$\widehat U_{\p(E)}$ and,  therefore, $\dist(z,E\setminus\widehat U_{\p(E)})\ge
\dist(A,E\setminus U)$.
Now we can adapt the proof of \cite[Proposition 2.4]{DinVen04} to
obtain our result.
$\square$ \end{proof}

Example \ref{ejemplo} below shows that if we drop off the assumption of boundedness, extensions to $\widehat U_{\p(E)}$ need not be of bounded type. However, it is possible to obtain extensions which are of bounded type ``around every point'' of $\widehat U_{\p(E)}$ in the following sense.

\begin{proposition}\label{extension de tipo acotada en cada punto}
Let $U$ be a balanced open set. For each $y\in\widehat U_{\p(E)}$, there exist a connected open subset $U_y$ of $\widehat U_{\p(E)}$ containing $U\cup \{y\}$ such that the extension to $\widehat U_{\p(E)}$ of any $f\in H_b(U)$  is of bounded type on $U_y$.
\end{proposition}

For the proof we will use the following two Lemmas, which are similar to some results in~\cite{din71(Cartan-Thullen)}. Although we state them for balanced open sets, they also hold for any $U$ for which polynomials are dense in $H_b(U)$. In the sequel, $\tilde f$ denotes the extension of $f\in H_b(U)$ to $\widehat U_{\p(E)}$.
\begin{lemma}
Let $U$ be a balanced open set, let $A$ a $U$-bounded set and let $y\in\widehat A_{\p(E)}$. Then for each $f\in H_b(U)$ we have
$$\|{d^k\tilde f(y)}\|\le\sup_{x\in A}\|{d^kf(x)}\|.$$
\end{lemma}
\begin{proof}
Given $\phi\in\p(^kE)'$, the function $\phi\circ{d^kf}$ is holomorphic and of bounded type on $U$. Its (unique) holomorphic extension to $\widehat U_{\p(E)}$ is given by $\phi\circ{d^k\tilde f}$. By Theorem~\ref{envoltura H_b balanceados} we have
$$\big|\phi\big({d^k\tilde f}(y)\big)\big| \le \|\phi\circ{d^kf}\|_A = \sup_{x\in A}\big|\phi\circ{d^kf(x)}\big|\le \|\phi\| \sup_{x\in A}\big|{d^kf(x)}\big|.$$ Since this is true for every $\phi\in\p(^kE)'$, the result follows.
$\square$ \end{proof}

\begin{lemma}
Let $U$ be a balanced open set, $A$ be a $U$-bounded set and let $3d=\dist_U(A)$. For   $z\in\widehat{A}_{\p(E)}$ we have  $B_E(z,d)\subset \widehat U_{\p(E)}$ and also $$\|\tilde f\|_{B_E(z,d)}\le\|f\|_{A+B_E(0,2d)}$$ for every  $f\in H_b(U)$.
\end{lemma}
\begin{proof} Let us write $C=A+B_E(0,2d)$.
For $x\in A$, the previous Lemma and Cauchy's inequalities imply
$$
\big\|\frac{d^k\tilde f(z)}{k!}\big\|\le \sup_{x\in A}\big\|\frac{d^kf(x)}{k!}\big\| \le \Big(\frac{1}{2d}\Big)^k \|f\|_{C}.
$$
Then, by the Cauchy-Hadamard formula, the Taylor series of
$\tilde f$ at $z$ converges in $B_E(y,2d)$. Also,
 for $\|x\|<d$ we obtain
$$
\sum_{k=0}^\infty \Big|\frac{d^k\tilde f(z)}{k!}(x)\Big|\le\sum_{k=0}^\infty \big\|\frac{d^k\tilde f(z)}{k!}\big\|\|x\|^k \le \|f\|_{C}\sum_{k=0}^\infty\Big(\frac1{2d}\Big)^kd^k=2\|f\|_{C}.
$$
This is true for every function in $H_b(U)$. In particular for each polynomial $P$ we have $\|P\|_{B_E(z,d)}\le2\|P\|_{C}$ and taking powers of $P$ we conclude that $\|P\|_{B_E(z,d)}\le\|P\|_{C}$. This implies that $B_E(z,d)\subset\widehat C_{\p(E)}\subset\widehat U_{\p(E)}$. Therefore, $$\sum_{k=0}^\infty \frac{d^k\tilde f(z)}{k!}(x)=\tilde f(z+x)$$ and  $\|\tilde f\|_{B_E(z,d)}\le\|f\|_{C}$.
$\square$ \end{proof}

Now we are ready to prove our proposition.

\begin{proof}(of Proposition \ref{extension de tipo acotada en cada punto})
Let $A$ be a $U$-bounded balanced set such that $y$ belongs to $\widehat{A}_{\p(E)}$. By  \cite[Lemma 1.4]{Vie07} the set  $\widehat{A}_{\p(E)}$ is balanced and hence it contains the segment $[0,y]$ joining $0$ and $y$.   By the previous Lemma, there exists a $U$-bounded set $C$ such that, for each $f\in H_b(U)$ and each $z\in[0,y]$, $$\|\tilde f\|_{B_E(z,d)}\le\|f\|_{C}<\infty,$$ where $3d=\dist_U(A)$. If we define $$U_y:=U\cup \Big(\bigcup_{z\in[0,y]}B_E(z,d)\Big), $$ then $y$ belongs to $U_y$ and $\tilde f$ is of bounded type on $U_y$.
$\square$ \end{proof}

Now we present an open balanced set $U\subset c_0$ and a function in $H_b(U)$ which cannot be extended to a holomorphic function of bounded type in $H_b(\widehat U_{\p(E)})$. By Theorem~\ref{envoltura H_b balanceados},  this answers for the negative the question made by Hirschowitz in \cite[Remarque~1.8]{Hir72}: is the extension of a function of bounded type to its $H_b$-envelope of holomorphy necessarily a function of bounded type?
Since by Theorem~\ref{envelopes coincide} the $H_b$-envelope is contained in the spectrum, this also shows that canonical extensions to the spectrum are not always of bounded type (see Section~\ref{seccion-bolalp} for more on this question).
This example is somehow inspired in  \cite[Example 7]{CarGarMae05}.

\begin{example}\label{ejemplo}
There exists an open balanced set $U\subset c_0$ and a function $g\in H_b(U)$ whose extension to  $\widehat U_{\p}$ is not of bounded type.
\end{example}
\begin{proof}
As a first step, let us see that it is enough to find an open balanced set $U\subset c_0$ satisfying:
\begin{enumerate}
            \item[(a)] For each $M>0$ there exists $k(M)\in \mathbb N$ such that $|x_{2n+1}|<3/4$ for every $x\in U$ with $\|x\|<M$ and for every $n\ge k(M)$.
            \item[(b)] The set $D:=\{e_{2n+1}:\, n\in\mathbb N\}$ is $\widehat U_{\p}$-bounded.
\end{enumerate}
For such a set $U$, we can define an entire function on $c_0$ by
$$ g(x):=\sum_{n\in\mathbb N} \Big(\frac{5}4 x_{2n+1} \Big)^n.$$
By property (b), $g$ is of bounded type in $U$. On the other hand,  $g(e_{2n+1}) \to \infty$ and this means, by (a), that $g$ is not of bounded type on  $\widehat U_{\p}$.

Now we show how to construct $U$. Define for $k,j\in\mathbb N$, $$p_{k,j}(x)= \big|kx_{2k+1}+x_{2j}\big| + \sup\{|x_{i}|:\, i\ne 2k+1,\; i\ne 2j\},$$
and $V_{k,j}=\{y\in c_0:\; p_{k,j}(y)<1/2\}$.
 Let us see that the balanced open set $$ U=\bigcup_{k,j\in\mathbb N}V_{k,j} + \frac{1}{4}B_{c_0}$$ satisfies properties (a) and (b) above.

To see (a), given $x\in U$ with $\|x\|\le M $, we can write $x=y+z$ with $4z\in B_{c_0}$ and  $y\in V_{k,j}$, for some $k,j\in\mathbb N$.  For any  $n\ne k$, we have $|y_{2n+1}|\le p_{k,j}(y)<1/2 $.  For $n=k$, we have  $$\frac12>p_{n,j}(y)\ge \big|ny_{2n+1}+y_{2j}\big|\ge
   n|y_{2n+1}|-M-\frac14$$ (since  $|y_{2j}|\le\|y\|\le M+1/4$), and then $|y_{2n+1}|<(M+1)/{n}<1/2$ whenever $n=k$ is greater than $2M+2$.
As a consequence, if $n\ge 2M+2$ we obtain $|x_{2n+1}|\le|y_{2n+1}|+|z_{2n+1}|<3/4$ (regardless of the relationship between $k$ and $n$).

To see (b), note that $p_{k,j}(e_{2k+1}-ke_{2j})=0$, and then $e_{2k+1}-ke_{2j}$ belongs to $V_{k,j}$ for every $k,j\in\mathbb N$. This implies that $e_{2k+1}-ke_{2j} + \frac14B_{c_0}$ is contained in $U$, and then the set $C_k:=\{e_{2k+1}-ke_{2j}+\frac18 B_{c_0}:\; j\in\mathbb N\}$ is $U$-bounded for every $k\in\mathbb N$.
Polynomials on $c_0$ are weakly sequentially continuous (see for example \cite[Proposition 1.59]{Din99})  and $e_{2j}$ converges weakly to $0$. Therefore, for each polynomial $P\in\p(c_0)$ and  for $\|z\|<1/8$  we have
$$
|P(e_{2k+1}+x)|\le \sup_{j\in\mathbb N}|P(e_{2k+1}-ke_{2j}+x)|\le \|P\|_{C_k}
$$
This means that $e_{2k+1}+x\in(\widehat{C_k})_{\p(c_0)}\subset\widehat U_{\p(c_0)}$ if $\|x\|<\frac18$ and (b) follows.
$\square$
\end{proof}

\smallskip

In view of the previous example, one may wonder if there exists some kind of envelope of $U$ (necessarily smaller than the $H_b$-envelope) to which the extensions of  functions in $H_b(U)$ are also of bounded type. The previous example together with Proposition~\ref{extension de tipo acotada en cada punto} suggest that we cannot expect the maximality property that envelopes are supposed to have. Let us see that even in the framework of Riemann domains we do not have this special envelope. The following result is widely known and follows from a straightforward connectedness argument (see, for example, \cite[Proposition 1.3]{Coe74}).
\begin{lemma}\label{igual o distinto}
Let $(X,p)$, $(Y,q)$ be connected Riemann domains spread over a Banach space $E$ and let $u,v:X\to Y$ be morphisms. Then either $u(x)=v(x)$ for every $x\in X$ or $u(x)\ne v(x)$ for every $x\in X$.
\end{lemma}

Now, suppose that  $(Y_0,\tau)$ is a Riemann domain which is, loosely speaking, the greatest among those with the following property: functions of $H_b(X)$ uniquely extend to functions in $H_b(Y_0)$. We put $V_n=\delta(X)\cup W_n^\circ$, where $W_n$ was defined in the proof of Theorem~\ref{envelopes coincide}. We have seen that the extension of every function in $H_b(X)$ to $V_n$ is of bounded type. The maximality of $Y_0$  gives morphisms $\nu_n:V_n\to Y_0$ such that $\tau=\nu_n\circ i_n$, where $i_n:X\to V_n$ is the inclusion for each $n$. For $m>n$, the application $\nu_m|_{V_n}:V_n\to Y_0$ is a morphism, and since $\nu_n|_U=\nu_m|_U=\tau$ we have $\nu_m|_{V_n}=\nu_n$ by Lemma~\ref{igual o distinto}.
Therefore, if $Y$ is the $H_b$-envelope of holomorphy of $X$,  we can define a morphism $\nu:Y\to Y_0$ by $\nu(x)=\nu_n(x)$ if $x\in V_n$.
On the other hand, it is clear that we have an $H_b$-extension morphism from $X$ to $Y_0$, so we have a morphism from $Y_0$ to $Y$ which is the inverse of $\nu$ by Lemma~\ref{igual o distinto}.
Therefore, such an $Y_0$ must coincide with $Y$, the already defined $H_b$-envelope of holomorphy of $X$. But we have seen that extensions to $Y$ are not necessarily of bounded type.

\section{Extending functions of bounded type to open subsets of $E''$}\label{seccion-bidual}

In this short section we consider extensions of functions of $H_b(U)$ to open sets in $E''$ containing $U$.
Let us start by defining the following variation of the set $\v U$ that also considers elements of the bidual:
$$
\v {U''}:=\{z\in E''\; :\; \textrm{there is some }\varphi\in M_b(U)\textrm{ such that }\varphi(f)=\overline f(z)\textrm{ for every }f\in H_b(E)\},
$$
where $\overline f$ denotes the Aron-Berner extension of $f$ \cite{AroBer78}.
Note that $\v{ U}=\v {U''}\cap E$ and $\v{ U''}\subset\pi(M_b(U))$ for every open set $U$. We also define for a $U$-bounded set $A$,
$$
\widehat A_{\p(E)}''=\{x''\in E'':\; |\overline f(x'')|\le \|f\|_A \;\; \textrm{for every } f\in\p(E)\}.
$$
We set
$$
\widehat U''_{\p(E)}:=\bigcup_{n\in\mathbb N}\Big(\widehat U_n\Big)''_{\p(E)}.
$$
\begin{remark}\label{techito prima prima y espectro}
We can prove as in Lemma \ref{propiedades de los techitos} that if $U$ is balanced then
\begin{equation*}
\widehat U''_{\p(E)}=\v {U''}.
\end{equation*}
Also, if $E$ is symmetrically regular, then $\widehat U''_{\p(E)}$ can be identified with the connected component of $M_b(U)$ which intersects $U$.
\end{remark}
Before we go on, let us make clear that we cannot expect $\widehat U_{\p(E)}''$ to be the largest open subset of $E''$ to which functions on $H_b(U)$ extend. For example, take a nonreflexive Banach space $E$ that is complemented in its bidual $E''$, say $E''=E\oplus M$. Denote by $\pi_E$ the projection to $E$. Then every function $f\in H_b(U)$ can be extended to $\tilde f\in H_b(U\times M)$ by $\tilde f=f\circ \pi_E$. On the other hand, the Hahn-Banach theorem shows that $\widehat U''_{\p(E)}\subset\overline{j_E(coe(U))}^{w^*}$ (where $j_E$ is the canonical inclusion of $E$ in its bidual and $coe(U)$ denotes  the absolutely convex hull of $U$). Thus, in general we can extend to sets which are larger than $\widehat U''_{\p(E)}$.
Things are different if we consider extensions that coincide locally with the Aron-Berner extension.

\begin{definition}\label{defi-abmorfismo}  If $W$ is an open subset of $E''$ containing $U$, a  continuous operator $e:H_b(U)\to H(W)$  will be called an \emph{$AB$-extension operator} if for some $x\in U$ (and hence for every $x\in U$)  there exists $r>0$ such that $e(f)$ coincides with the Aron-Berner extension of $f$ on $B_{E''}(x,r)$.
\end{definition}

The following proposition shows that  $\widehat U''_{\p(E)}$ can be seen as an $AB$-envelope of $U$, at least for balanced open sets. This result can be seen as  analogous to Theorems \ref{envoltura H_b balanceados} and \ref{extension acotados}, modulo Aron-Berner extensions, and the proof is similar, using Remark~\ref{techito prima prima y espectro}.

\begin{proposition}\label{extension E segunda mas}
Let $U$ be an  balanced open subset of a symmetrically regular Banach space.  \newline
($a$) There exists an {$AB$-extension operator} from   $H_b(U)$ to $H(\widehat U''_{\p(E)})$. The extension of each function in $H_b(U)$ is bounded on the sets $\widehat A''_{\p(E)}$, for every $U$-bounded set $A$. Also, for $z\in \widehat A''_{\p(E)}$, we have $\delta_z\prec A$. \newline
($b$) If  in addition $U$ is bounded, extensions belong to $H_b(\widehat U''_{\p(E)})$. Also, $\widehat U''_{\p(E)}$ is a $H_b$-domain of holomorphy.
\end{proposition}

We remark that, if there exists an {$AB$-extension operator} $e:H_b(U)\to H(W)$ for some $W\subset E''$, then $W$ must be a subset of $\widehat U''_{\p(E)}$. Indeed, if $z\in W\setminus \widehat U''_{\p(E)}$,  just as in Theorem~\ref{envelopes coincide} we can choose functions $f_n\in H_b(E)$ such that $|\overline f_n(z)|>1$ and  $f_n\to 0$ in $H_b(U)$. Thus
$e(f_n)\to0$ in $H(W)$, which is impossible since $|e(f_n)(z)|=|\overline f_n(z)|>1$ for all $n$.

Davie and Gamelin \cite{DavGam89} showed that the Aron-Berner extension is an isometry from $H^{\infty}(B_E)$ to $H^{\infty}(B_{E''})$. Later, it was shown in \cite[Theorem 1.3]{GalGarMae93} that if $U$ is convex and balanced then the Aron-Berner extension is isometric isomorphism from $H^\infty(U)$ to $H^\infty(int(\overline{U}^{w^*}))$, where $int(\overline{U}^{w^*})$ means the norm-interior of the weak-star closure of $U$ in $E''$. Moreover, Theorem 1.5 in \cite{GalGarMae93} asserts that there exists an $AB$-extension morphism (in the sense of Definition~\ref{defi-abmorfismo}) from $U$ to  $int(\overline{U}^{w^*})$. As a consequence, $int(\overline{U}^{w^*})$ is contained in $\widehat U''_{\p(E)}$.
The reverse inclusion is an easy consequence of the Hahn-Banach theorem, so we have that $\widehat U''_{\p(E)}=int(\overline{U}^{w^*})$ for any convex and balanced $U$.
\medskip

\section{Density of finite type polynomials}\label{seccion-density}

In several complex variables, the  holomorphic convexity of $U$, or $U$ being a domain of holomorphy, is equivalent to $M_b(U)=\delta(U)$. In our infinite dimensional setting this is not the case unless $E$ has very particular properties. These particular properties arise rather naturally. If $E$ is not reflexive, there are always elements of the bidual in the spectrum, so the equality $M_b(U)=\delta(U)$ cannot hold. On the other hand, if there are polynomials on $E$ that are not weakly continuous on bounded sets, there is much more than evaluations in the spectrum \cite{AroColGam91,AroGalGarMae96}, and so $M_b(U)=\delta(U)$ is impossible regardless of the reflexivity of $E$. We will formalize this below, refining some results of \cite{Vie07,Muj01}.

Recall that a \emph{Tsirelson-like space} is a reflexive Banach space on which every polynomial is approximable (i.e, limit in norm of finite type polynomials).  For Tsirelson-like spaces, it is clear that $M_b(E)$ identifies with $E$.
In \cite[Theorem 2.1]{Vie07}, Vieira proved that if $U$ is a balanced $H_b(U)$-convex subset of a Tsirelson-like space, then $M_b(U)=\delta(U)$.
This can be seen in the following way. If every polynomial on $E$ is approximable and $U$ is balanced, then finite type polynomials are dense in $H_b(U)$. This implies that $\pi$ is injective and $M_b(U)$ is contained in $\delta(E'')$. By Remark~\ref{techito prima prima y espectro}, $M_b(U)$ must coincide with $\delta(\widehat U_{\p(E)}'')$. If $E$ is also reflexive, then we have $M_b(U)=\delta(\widehat U_{\p(E)})$. This last set coincides with $U$ whenever $U$ is $H_b(E)$-convex and, in particular, when it is $H_b(U)$-convex. This shows \cite[Theorem 2.1]{Vie07}.
As a consequence of our results, we can also show a kind of converse of this result.

\begin{theorem}\label{tsirelson_reciproco}
Let $U$ be a balanced open subset of a Banach space $E$ with the approximation property. Then,  $M_b(U)=\delta(U)$ if and only if $U$ is $H_b(U)$-convex and $E$ is a Tsirelson-like space.
\end{theorem}
\begin{proof}
If  $M_b(U)=\delta(U)$, by Lemma~\ref{propiedades de los techitos} we have
$$
\delta(U)\subset\delta(\widehat U_{\p(E)})\subset M_b(U) =\delta(U).
$$
Therefore, $\widehat U_{\p(E)}=U$ and $U$ is $\p(E)$-convex. Since $U$ is  balanced, then it is $H_b(U)$-convex~\cite[Proposition 1.5]{Vie07}.  With the same proof of  Theorem 1.2  in \cite{Muj01}, we can see that the equality $M_b(U)=\delta(U)$ implies that  $E$ is reflexive and every polynomial on $E$ is weakly continuous on bounded sets. Now, since $E$ is reflexive and has the approximation property, so does $E'$ and polynomials on $E$ are then approximable \cite[Proposition 2.7]{AroPro80}.

The converse is the already mentioned Theorem 2.1 of \cite{Vie07}.
$\square$ \end{proof}

\smallskip
As the previous theorem states, the equality $M_b(U)=\delta(U)$ is hard to achieve for domains in an arbitrary Banach space $E$. This is mainly because  $M_b(U)$ cannot be, in general, identified with a subset of $E$. But we know that $M_b(U)$ can be projected on $E''$ via $\pi$, so a natural question is the following: suppose that $U$ is  $H_b(U)$- or $H_b(E)$-convex and $E$ is reflexive. Is it true that  $\pi(M_b(U))=U$? And if we drop off the reflexivity assumption, can we obtain something like $\pi(M_b(U))=\v{ U''}$ instead?
The answer relies on the density of finite type polynomials. Namely, we will see that if there is a non approximable polynomial on $E$, there are proper $H_b(E)$-convex subsets $U$ of $E$ for which $\pi(M_b(U))$ is larger than $\v{ U''}$, since it contains the whole space $E$. In particular, if $E$ is reflexive with the approximation property but not Tsirelson-like, there are subsets $U\subsetneq E$ that are $H_b(E)$-convex satisfying $\pi(M_b(U))=E$. First we state and prove the following easy result.

\begin{lemma}\label{distanciapositiva} Let $C$ and $D$ be subsets of a Banach space $E$, one of them bounded. If there exist a homogeneous polynomial $P$ on $E$ and $\varepsilon >0$ such that $|P(x)-P(y)|>\varepsilon$ for all $x\in C$ and $y\in D$, then the distance between $C$ and $D$ is strictly positive.
\end{lemma}
\begin{proof} Let us assume that $C$ is bounded and let $R$ be the radius of a ball containing $C$. If $\|y\|> R +1$, then clearly $\|x-y\|>1$ for all $x\in C$. On the other hand, the polynomial $P$ is uniformly continuous on the closed ball $\bar B(0,R+1)$. Therefore, since $|P(x)-P(y)|>\varepsilon$ for all $x\in C, y\in D$, there must exist $\delta>0$ such that $\|x-y\|>\delta$ for $x\in C$, $y\in D\cap \bar B(0,R+1)$. This completes the proof. $\square$ \end{proof}

If $C$ is a bounded subset of $U$ and $D$ is $E\setminus U$, then the existence of $P$ and $\varepsilon$ as in the lemma ensures that $C$ is $U$-bounded. We remark that the conclusion of the previous lemma does not hold if  $C$ and $D$ are unbounded. Consider for example in $\mathbb C^2$ the sets $C= \{(x,y):xy\ge1\}$ and $D=\{(x,y): y=0\}$, together with the polynomial $P(x,y)=xy$.

\begin{proposition}\label{tl}
Suppose $E'$ has the approximation property. The following conditions are equivalent: \newline
(i) finite type polynomials are dense in $\p(E)$; \newline
(ii) for every open subset $U$ of $E$ we have $\v{ U''}=\pi(M_b(U))$;\newline
(iii) for every open $H_b(E)$-convex subset $U$ of $E$ we have $\v{ U''}=\pi(M_b(U))$.

\noindent If the conditions do not hold, then there exists a proper subset $U$ of $E$ which is $H_b(E)$-convex but which satisfies $\pi(M_b(U))\supset E$.
\end{proposition}
\begin{proof}
($i$) $\Rightarrow$ ($ii$)
Let $z=\pi(\varphi)$ for some  $\varphi\in M_b(U)$. Since finite type polynomials are dense in $H_b(E)$ and $H_b(E)$ continuously embeds in $H_b(U)$, we have $\varphi(f)=\overline f(z)$ for every $f\in H_b(E)$, where $\overline f$ denotes the Aron-Berner extension of $f$. As a consequence, $z\in \v{U''}$ by the very definition of this set. The reverse inclusion is easy.

Clearly,  ($ii$) implies ($iii$).

($iii$) $\Rightarrow$ ($i$) If finite type polynomials are not dense in $\p(E)$, since $E'$ has the approximation property, there must exist a $K$-homogeneous polynomial $P$ which is not weakly continuous on bounded sets \cite[Proposition 2.7]{AroPro80}. Define the set
$$
U=\{x\in E:\; Re(P(x))>\frac12\}.
$$
Let us see that $U$ is $H_b(E)$-convex, then that $\v{U''}\cap E=U$,  and finally that  $\pi(M_b(U))$ contains $E$.

We consider a  fundamental sequence of $U$-bounded sets  $(U_n)_n$ as in (\ref{typical fundamental sequence}).
For  fixed $n\in \mathbb N$ and $x\in U_n$ we set  $$\alpha =\left(\frac{1}{2Re(P(x))}\right)^\frac{1}{K}<1.$$
We have $Re(P(\alpha x))=\alpha^K Re(P(x))=\frac12$, so  $\alpha x$ does not belong to $U$. This means that  $\|x-\alpha x\|\ge \frac{1}{n}$, from which we get $$1-\left(\frac{1}{2Re(P(x))}\right)^\frac{1}{K}\ge\frac{1}{n\|x\|}\ge\frac{1}{n^2},$$and then \begin{equation}\label{eldeltan}\eta_n:= \frac12\left(\frac{n^2}{n^2-1}\right)^K \le Re(P(x)).\end{equation} We have shown that $Re(P(x))\ge \eta_n>1/2$ for all $x\in U_n$.
As a consequence, $|e^{-P(x)}|\le e^{-\eta_n}$ for all $x\in U_n$ and the same must hold for any $x\in \widehat {(U_n)}_{H_b(E)}$. Therefore, $Re(P(x))-Re(P(y))>\eta_n-1/2>0$ for all $x\in \widehat {(U_n)}_{H_b(E)}$ and all $y\in E\setminus U$. Since $\widehat {(U_n)}_{H_b(E)}$ is clearly bounded, it turns out to be $U$-bounded by Lemma~\ref{distanciapositiva}. This holds for any $n$ and we conclude that $U$ is $H_b(E)$-convex.
By the definitions of the sets $\v{U''}$ and $\v U$ we have $\v{U''}\cap E=\v U$, and this last set coincides with $U$, since $U$ is $H_b(E)$-convex (see Lemma~\ref{propiedades de los techitos}).

The final step is to show that $E\subset\pi(M_b(U))$. Take a weakly null  net $\{x_i\}_{i\in I}\subset S_E$ such that $P(x_i)$ is a real number greater than 1 for every $i\in I$. For $x\in E$, we can find $\lambda>0$ such that the bounded set $\{x+\lambda x_i\}_{i\in I}$ is actually $U$-bounded. Indeed, for each $i\in I$, the mapping $\lambda\mapsto P(x+\lambda x_i)$ is a polynomial of degree $K$ whose leading term is $P(x_i)\lambda^K$. Since $P(x_i)>1$ for every $i$, we can choose $\lambda$ large enough to get $Re(P(x+\lambda x_i))\ge 1$ for all $i\in I$. Fixed such $\lambda$, for every  $y\in E\setminus U$ we have $|P(x+\lambda x_i)-P(y)|\ge |Re(P(x+\lambda x_i))-Re(P(y))| \ge 1/2$. Lemma~\ref{distanciapositiva} ensures that $\{x+\lambda x_i\}_{i\in I}$ is $U$-bounded.
Then $\{x+\lambda x_i\}_{i\in I}$ is contained in $U_N$ for some $N>0$. Since $\{x_i\}_{i\in I}$ is weakly null, we have $x\in \overline{U_N}^{w^*}$ and, by \cite[Proposition 18]{CarGarMae05}, we obtain $x\in\pi(M_b(U))$. This holds for any $x\in E$, so we conclude that $E\subset \pi(M_b(U))$.
$\square$ \end{proof}

\bigskip

      We end this section with a Banach-Stone type result. First we fix some notation.
      For an open set $U\subset E$ and a family $\mathcal A=(A_k)_k$ of subsets of $U$ with $\bigcup_k A_k=U$, we define as in~\cite{Muj86,MujNac92}
      $$
      \mathcal H^\infty\big( \mathcal A \big)=\{f\in H(U):\; \|f\|_{A_k}<\infty \textrm{ for every }k\}.
      $$
      This is a Fr\'echet algebra with the topology of uniform convergence on the $A_k$'s. If $\mathcal A$ is a fundamental system of $U$-bounded sets, then we have simply $\mathcal H^\infty\big(\mathcal A\big)=H_b(U)$. Note that, if $U$ is a balanced subset of a symmetrically regular Banach space $E$, Proposition~\ref{extension E segunda mas} states that every function $f\in H_b(U)$ can be extended to a function $\overline f\in\mathcal H^\infty\big(\mathcal U\big)$, where $\mathcal U= ( (\widehat{U_k})''_{\p(E)})_k$.

      If $V\subset F$ is an open set and we have a family   $\mathcal B=(B_j)_j$  of subsets of $V$ such that $\bigcup_j B_j=V$, we define the Fr\'echet algebra
      $$
      \mathcal H^\infty\big( \mathcal B,\mathcal A\big)  =  \{g\in H(V,U):\; \textrm{ for each } j,\  g(B_j) \textrm{ is contained in some } A_{k} \}.
      $$
      If $\mathcal A$ and $\mathcal B$ are fundamental systems of $U$-bounded sets and $V$-bounded sets respectively, then $\mathcal H^\infty\big( \mathcal B, \mathcal A\big)$ is the algebra $H_b(V,U)$ of holomorphic functions of bounded type from $V$ to $U$.

We will say that a function $g$ defined on a dual Banach space (with values in some topological space) is \textit{locally $w^*$-continuous} at some point if there exists a (norm) neighborhood of it such that the restriction of the function to this neighborhood is $w^*$-continuous. A function is locally $w^*$-continuous on an open set if it is locally $w^*$-continuous at each point of the set. Also, a function between dual spaces is \textit{locally $w^*$-$w^*$-continuous} if it is locally $w^*$-continuous when range space is endowed with the weak-star topology.

      In \cite[Theorem 3.1]{Vie07} the following is proved: if $E$ and $F$ are reflexive Banach spaces, one of them Tsirelson-like, and $U\subset E$, $V\subset F$ are open balanced and $\p(E)$-convex, then
      the following conditions are equivalent:
      \begin{itemize}
      \item[($i$)] There exists a bijective mapping $g:V\to U$ such that $g\in H_b(V,U)$  and $g^{-1}\in H_b(U, V)$.
      \item[($ii$)] The algebras $H_b(U)$ and $H_b(V)$ are topologically isomorphic.
      \end{itemize}
      When the conditions are satisfied, then $E$ and $F$ are isomorphic Banach spaces.
      In \cite[Corollary 22]{CarGarMae05} a similar result was proved for convex balanced open sets when every polynomial on either $E''$ or $F''$ is approximable. In which case it follows that $E'$ and $F'$ are isomorphic.

      We slightly improve these results with the Banach-Stone type result stated in Theorem~\ref{Banach-Stone}. In the rest of this section,   $(U_k)_k$ and $(V_j)_j$ will stand for fundamental systems of $U$-bounded and $V$-bounded sets respectively, and we define $\mathcal U= ( (\widehat{U_k})''_{\p(E)})_k$ and $\mathcal V= ( (\widehat{V_k})''_{\p(F)})_k$. In the following Theorem and Lemmas,  the bar indicates the $AB$-extension of a function to the corresponding open subset of the bidual.

      \begin{theorem}\label{Banach-Stone}
      Let $E,F$ be Banach spaces, $V\subset F$, $U\subset E$ open balanced subsets and suppose that every polynomial on $E''$ is approximable.
      If $\phi:H_b(U)\to H_b(V)$ is a Fr\'echet algebra isomorphism then there exists a biholomorphic function $g:\widehat V''_{\p(F)}\to \widehat U''_{\p(E)}$, with $g\in\mathcal H^\infty\big( \mathcal V,\mathcal U\big)$ and $g^{-1}\in\mathcal H^\infty\big(\mathcal U,\mathcal V\big)$, both locally $w^*$-$w^*$-continuous, such that $\overline{\phi f}=\overline f\circ g$ for every $f\in H_b(U)$.
      Conversely, if $g$ is such a function then the operator $\phi:H_b(U)\to H_b(V)$ given by $\phi f=\overline f\circ g|_V$ is a Fr\'echet algebra isomorphism.

      In any of these cases, $E'$ and $F'$ are isomorphic Banach spaces.
      \end{theorem}

      To prove this Theorem we will need some preliminary results.
      Let $V\subset F$ be a balanced subset. Aron and Berner's result \cite[Corollary 2.1]{AroBer78}, together with the isometry proved by Davie and Gamelin in~\cite{DavGam89}, show that there is an $AB$-extension operator $f\mapsto \bar f$ from $H_b(V)$ to $H(W)$, where $W$ is the subset of $E''$ given by
\begin{equation}\label{el_W}
W=\bigcup_{x\in U} B_{E''}(x,\dist(x,F\setminus V)).
\end{equation}
Moreover, the application $\delta_{y''}(f):=\overline f(y'')$ defines a continuous homomorphism on $H_b(V)$ for each $y''\in W$.

      \begin{lemma}\label{la g es holo}
      Let $E,F$ be Banach spaces, $V\subset F$ an open balanced subset and $U\subset E$ open. Suppose that $\phi:H_b(U)\to H_b(V)$ is a continuous, linear and multiplicative operator. Then

      $(a)$ the mapping $g:W\to E''$, defined by $g(y'')=\pi(\delta_{y''}\circ \phi)$ is holomorphic (where $W$ is the set defined in (\ref{el_W}));

      $(b)$ if $F$ is symmetrically regular, then the mapping $g:\widehat V''_{\p(F)}\to E''$, defined by $g(y'')=\pi(\delta_{y''}\circ \phi)$ is holomorphic.
      \end{lemma}
      \begin{proof}
      $(a)$ Denote by $\theta_\phi:M_b(V)\to M_b(U)$ the transpose of $\phi$ restricted to the spectra. Then $g$ is just the composition $W \overset{\delta}{\longrightarrow} M_b(V) \overset{\theta_\phi}{\longrightarrow} M_b(U) \overset{\pi}{\longrightarrow} E''$, which is well defined by the comments above. If we take $y''\in W$ and $x'\in E'$, then $g(y'')(x')=\delta_{y''}(\phi x')=\overline{\phi x'}(y'')$. Thus $g$ is weak*-holomorphic on $W$ and therefore holomorphic (see for example \cite[Exercise 8D]{Muj86}).

      $(b)$ By Proposition \ref{extension E segunda mas} we can define $g$ on $\widehat V''_{\p(F)}$ and proceed as in the proof of $(a)$.
      $\square$ \end{proof}

      \begin{lemma}\label{operador de composicion}
      Suppose that every polynomial on $E$ is approximable and let $F$ be symmetrically regular. If $V\subset F$, $U\subset E$ are open balanced subsets and $\phi:H_b(U)\to H_b(V)$ is a continuous linear operator, then
      $\phi$ is multiplicative if and only if there exists a holomorphic function $g:\widehat V''_{\p(F)}\to \widehat U''_{\p(E)}$ belonging to $\mathcal H^\infty\big(\mathcal V,\mathcal U\big)$ such that $\overline{\phi f}=\overline f\circ g$ for every $f\in H_b(U)$.
      \end{lemma}
      \begin{proof} First we note that, if every polynomial on $E$ is approximable, then $E$ must  be symmetrically regular.
      Suppose that $\phi$ is multiplicative and let $g$ be the mapping defined by Lemma \ref{la g es holo} $(b)$. Since every polynomial on $E$ is approximable, the spectrum $M_b(U)$ can be identified with $\widehat U''_{\p(E)}$, thus $g$ maps $\widehat V''_{\p(F)}$ inside $\widehat U''_{\p(E)}$. Also, the definition of $g$ ensures that $\overline f(g(y''))=\overline{\phi f}(y'')$ for $f\in H_b(U)$ and $y''\in \widehat V''_{\p(F)}$.
      It remains to prove that $g$ belongs to $\mathcal H^\infty\big(\mathcal V,\mathcal U\big)$. Suppose by contradiction that, for some $n_0\in\mathbb N$, $g((\widehat{V_{n_0}})''_{\p(F)})$ is not contained in any of the $(\widehat{U_k})''_{\p(E)}$'s. So we can choose, for each $k$, an element  $y_k''$ in $ (\widehat{V_{n_0}})''$ such that $g(y_k'')$ is not in $(\widehat{U_k})''_{\p(E)}$. This means that there are polynomials $P_k\in\p(E)$ such that $\|P_k\|_{U_k}<1/{2^k}$ and $$|\overline P_k(g(y_k''))|>k+\sum_{j=1}^{k-1}|\overline P_j(g(y_k''))|.$$ Let $f=\sum P_k\in H_b(U)$. By Proposition~\ref{extension E segunda mas}, the function  $\overline{\phi f}$ belongs to $\mathcal H^\infty\big(\mathcal V\big)$. But this is a contradiction, since for any $k$ we have $\|\overline{\phi f}\|_{(\widehat{V_{n_0}})''_{\p(F)}}\ge |\overline{\phi f}(y_k'')|=|\overline f(g(y_k''))| >k-1$.

      The converse is immediate.
      $\square$ \end{proof}

      \bigskip Now we are ready to prove our Banach-Stone type result.

      \begin{proof}(of Theorem \ref{Banach-Stone})
      Suppose that $\phi$ is an isomorphism and let $g:W\to E''$ be holomorphic function given in Lemma~\ref{la g es holo} $(a)$. Our hypothesis imply that $E$ is symmetrically regular, so we can also consider  $h:\widehat U''_{\p(E)}\to F''$ the holomorphic map obtained from the homomorphism $\phi^{-1}$, using Lemma \ref{la g es holo} $(b)$. Then $h\circ g$ is the composition
      $$
	W\overset{\delta}{\to} M_b(V) \overset{\theta_\phi}{\to} M_b(U) \overset{\pi}{\to}  \widehat U''_{\p(E)} \overset{\delta}{\to} M_b(U) \overset{\theta_{\phi^{-1}}}{\to} M_b(V)\overset{\pi}{\to} F''.
      $$
      Since  $M_b(U)=\delta(\widehat U''_{\p(E)})$, we have $h\circ g=id_{W}$ and by differentiation,  $dh(g(0))\circ dg(0)=id_{F''}$. This means that $F''$ is isomorphic to a complemented subspace of $E''$, which implies that every polynomial on $F''$ is approximable and, as a consequence, that $F$ is symmetrically regular. Thus we can use Lemma \ref{operador de composicion} to define $g$ on $\widehat V''_{\p(F)}$,  and we have  $h\circ g=id_{\widehat V''_{\p(F)}}$. Since polynomials on $F'' $ are approximable, so are polynomials on $F$, and then we have  $M_b(V)=\delta(\widehat V''_{\p(F)})$. Thus $g\circ h=id_{\widehat U''_{\p(E)}}$ which means that $h=g^{-1}$. By Lemma \ref{operador de composicion}, both $g$ and $g^{-1}$ belong to the corresponding $\mathcal H^\infty$ algebras. They are also locally $w^*$-$w^*$ continuous. Indeed, for every  $x'\in E'$, $x'\circ g=\overline{\phi x'}$ is locally $w^*$-continuous since it is locally an Aron-Berner extension, and therefore $g$ is locally $w^*$-$w^*$ continuous on $V$. For $g^{-1}$ we can proceed analogously.

      Conversely, suppose that $g$ is as above. If we define $\phi f=\overline f\circ g|_V$ for $f\in H_b(U)$ and $\psi h=\overline h\circ g^{-1}|_U$ for $h\in H_b(V)$,  then clearly $\phi:H_b(U)\to H_b(V)$ and $\psi:H_b(V)\to H_b(U)$ are continuous and multiplicative operators. Let us see that $\psi=\phi^{-1}$. For  $f\in H_b(U)$ we have
      \begin{equation}\label{eqB-S}
      \psi\circ\phi f=\psi(\overline f\circ g|_V)= \overline{\overline f\circ g|_V}\circ g^{-1}|_U.
      \end{equation}
      Since every polynomial on $E$ is approximable, Aron-Berner extensions coincide locally with extensions by $w^*$-continuity and density. Therefore, $\overline f$ is locally $w^*$-continuous, and so is $\overline f\circ g$. For $z\in \widehat V''_{\p(F)}$, we can apply Lemma 2.1 of \cite{AroColGam95} to the restriction of $\overline f\circ g$ to a suitable ball, to obtain that ${d^k(\overline f\circ g)(z)}$ is a $w^*$-continuous polynomial for every $k$. By \cite[Theorem 2]{Zal90} we can conclude that $ {d^k(\overline f\circ g)(z)}$ is in the image of the Aron-Berner extension, and thus $\overline{\overline f\circ g|_V}=\overline f\circ g$. Substituting in (\ref{eqB-S}) we obtain $\psi\circ\phi (f)=f$. Similarly, we can show that $\phi\circ\psi ( h)=h$ for $h\in H_b(V)$.

      It remains to prove that $E'$ and $F'$ are isomorphic.
      As above, differentiating $g\circ g^{-1}$ at $0$ we obtain that $E''$ and $F''$ are isomorphic.
      An application of \cite[Lemma 2.1]{AroColGam95} to the restriction of $y''\mapsto
      g(y'')(x')$ to a suitable ball, gives that the
      differential of $g$ at any point is $w^*$-$w^*$-continuous. The same holds for the differential of $g^{-1}$ at any point. Therefore, the isomorphism between $E''$ and $F''$ is actually the transpose of an isomorphism between $F'$ and $E'$.
      $\square$ \end{proof}

When one of the spaces is Tsirelson-like and the open sets are bounded, we can say something more. In the following, the tilde denotes the extension of a function to the $\p(E)$-hull as in Theorem~\ref{envoltura H_b balanceados}.
      \begin{corollary}
      Let $E,F$ be Banach spaces, one of them Tsirelson-like, and let $V\subset F$, $U\subset E$ be open, balanced and bounded subsets. Then $\phi:H_b(U)\to H_b(V)$ is a Fr\'echet algebra isomorphism if and only if there exists a biholomorphic function $g\in H_b(\widehat V_{\p(F)},\widehat U_{\p(E)})$ such that $g^{-1}\in H_b(\widehat U_{\p(E)},\widehat V_{\p(F)})$ satisfying  $\widetilde{\phi f}=\tilde f\circ g$ for every $f\in H_b(U)$.
      In which case $E$ and $F$ are isomorphic Banach spaces.
      \end{corollary}

      The Tsirelson-James space $T_J^*$ is not reflexive (it is not a Tsirelson-like space) but satisfies the conditions of Theorem~\ref{Banach-Stone} by \cite[Lemma 19]{DimGalMaeZal04}.

\section{On the Spectrum of $H_{b}(U)$}\label{seccion-bolalp}

A consequence of Example \ref{ejemplo} is that the canonical extension of a function in $H_b(U)$ to the spectrum is not necessarily of bounded type. An application of the Closed Graph Theorem gives an equivalent condition for these extensions to be of bounded type. Suppose $U$ is an open subset of a symmetrically regular Banach space. Then the canonical extension of every  function in $H_{b}(U)$ to $M_{b}(U)$ is of bounded type if and only if given any $M_{b}(U)$-bounded set $B$ there
exists a $U$-bounded set $D$ such that  $\varphi\prec D$ for every $\varphi\in B$. See
\cite[Proposition 2.5]{DinVen04} for a related result.

In \cite{AroGalGarMae96} the following inequality was implicitly shown:
\begin{equation}\label{desigualdad dist}\sup\big\{\dist(A,U^{c}):\,A\subset U, \,\varphi\prec A\big\} \le \dist_{M_{b}(U)}(\varphi),\end{equation}where the meaning of $\dist_{M_{b}(U)}$ is given in~(\ref{defi-distancia}).
If for some $U$ we have equality or at least a reverse inequality with some constant, then extensions to $M_b(U)$ would be of bounded type, as a consequence of the above comments.  We do not know of many examples in which extensions to $M_b(U)$ are of bounded type. A first example is to take $U$ as $c_0$ or any space $E$ for which finite type polynomials are dense in $\p(E)$ (the original Tsireslon space $T^*$ is another example). In this case we have $M_b(E)=E''$ and the extension to the spectrum is the Aron-Berner extension, which is of bounded type.
More generally, if $U$  is any bounded and balanced open subset of such a space $E$, we have seen at the beginning of Section~\ref{seccion-density} that $M_b(U)$ coincides with
$\widehat U_{\p(E)}''$. Therefore, Proposition~\ref{extension E segunda mas} (b) ensures that extensions to $M_b(U)$ are of bounded type.
On the other hand, if $E$ is any symmetrically regular Banach space, it was shown in \cite[Proposition 6.30]{Din99} that the extension to the spectrum is of bounded type on each connected component of it. We now show that, in spite of this fact, these extensions need not be of bounded type on the whole spectrum.

\begin{proposition}\label{la extension al espectro}
Let $E$ be a symmetrically regular Banach space and suppose there exists a continuous $n$-homogeneous polynomial $P$ on $E$ which is not weakly continuous on bounded sets. Then the canonical extension of $P$ to the spectrum $M_b(E)$ is not of bounded type.
\end{proposition}
\begin{proof}
By \cite[Corollary 2]{BoyRya98}, the restriction of $P$ to some ball is not weakly continuous at 0. Then, we can take  a weakly null bounded net $(x_\alpha)_{\alpha\in\Delta}$  and $\varepsilon>0$ such that $|P(x_\alpha)|>\varepsilon$ for every $\alpha\in\Delta$.
For each $k\in\mathbb N$, let $\varphi_k\in M_b(E)$ be an adherent point of of the net $\{\delta_{kx_\alpha}\}_{\alpha\in\Delta}$. Since $(k x_\alpha)_{\alpha\in\Delta}$ is weakly null, we have $\varphi_k(x')=0$ for every $x'\in E'$,  and thus $\pi(\varphi_k)=0$ for every $k$. This implies that the set $\{\varphi_k:\,k\in\mathbb N\}$ is $M_b(E)$-bounded. However, $|\tilde P(\varphi_k)|=|\varphi_k(P)|>k^n \varepsilon$, and therefore $\tilde P$ is not of bounded type on $M_b(E)$.
$\square$ \end{proof}

Most classical Banach spaces admit polynomials satisfying the hypotheses of Proposition~\ref{la extension al espectro}, exceptions being $c_0$ and the original Tsirelson space $T^*$. For $E=\ell_p$ we simply take $P(x)=\sum x_j^n$ for some $n>p$.
Recall that the open set in Example \ref{ejemplo} was neither bounded nor convex, so one might ask if for the unit ball of a symmetrically regular Banach space things are easier. We do not know if in this case extensions to the spectrum are of bounded type, but we can answer for the negative the question on the reverse inequality in (\ref{desigualdad dist}): fixed $1<p<\infty$, there cannot be a constant $c>0$ such that $\sup\left\{\dist(A,B_{\ell_p}^{c}):\,\varphi\prec A\right\} \ge  c\ \dist_{M_{b}(B_{\ell_p})}(\varphi)$ for every $\varphi\in M_b(B_{\ell_p})$.

For the following proposition, we recall once again some facts on the analytic structure of $M_b(U)$ from~\cite{AroGalGarMae96},  this time for the specific case $U=B_{\ell_p}$, with $1<p\le\infty$. Fix $\varphi\in M_b(B_{\ell_p})$. If  $\varphi \prec rB_{\ell_p}$ for some $0<r<1$, then the homomorphism $\varphi^z\in M_b(B_{\ell_p})$ given by   $$\varphi^z(f)=\sum_{n=0}^\infty \varphi\left(y\mapsto \frac{ d^nf(y)}{n!}(z)\right)  $$
is defined for any $z\in\ell_p$  with $\|z\|<\frac1{1-r}$ (for the case $p=\infty$, for any $z\in\ell''_\infty$).
In fact, inequality~(\ref{desigualdad dist}) is a consequence of this: since $\varphi^z(f)$ is defined whenever $\|z\|<\frac1{1-r}$, we have $\dist_{M_{b}(B_{\ell_p})}(\varphi)\ge \frac1{1-r}$, which is precisely the distance from  $rB_{\ell_p}$ to $\ell_p\setminus B_{\ell_p}$.
In the sequel, for $\varphi \in M_b(B_{\ell_p})$ we define $S(\varphi)$, the sheet of $\varphi$, as the connected component of  $M_b(B_{\ell_p})$ that contains $\varphi$.

\begin{proposition}\label{no reverse ineq}
If $1<p<\infty$, then for each $\delta>0$ there exists $\varphi\in M_b(B_{\ell_p})$ such that
$$\sup \left\{\dist(A,B_{\ell_p}^{c}):\,\varphi\prec A\right\}< \delta\  {\dist_{M_{b}(B_{\ell_p})}(\varphi)}.$$
In other words, there is no reverse inequality in (\ref{desigualdad dist}).
\end{proposition}
\begin{proof} For each $0<s<1$, let $\varphi_s\in M_b(B_{\ell_p})$ be an adherent point of the set $\{\delta_{se_n}\}_n$.
Take a natural number $m>p$ and define $$g_N(x)=\Big(\sum_k \Big(\frac{x_k}s\Big)^m\Big)^N.$$ Then $\varphi_s(g_N)=1$ for all $N\in \mathbb N$. If $0<r<s$, we have $$\|g_N\|_{rB_{\ell_p}}\le\Big(\frac r s \Big)^{mN}\rightarrow 0$$ as $N\to\infty$. Thus $\varphi_s\nprec rB_{\ell_p}$ for $r<s$, and since we clearly have $\varphi_s\prec s B_{\ell_p}$, we obtain $\sup\left\{\dist(A,B_{\ell_p}^{c}):\,\varphi_s\prec A\right\}=1-s$.

On the other hand, take $z\in \ell_p$ with $\|z\|^p<1-s^p$. We can choose $\varepsilon>0$ and $n_0\in\mathbb N$ such that the set $A=\{\lambda z+s {e_n}:\, \, |\lambda|=1+\varepsilon,\, n\ge n_0 \}$ is $B_{\ell_p}$-bounded. By  Cauchy's inequalities we have
$$
\left|\frac{d^kf(se_n)}{k!}(z)\right|\le \frac1{(1+\varepsilon)^k}\sup_{|\lambda|=1+\varepsilon}|f({s e_n}+\lambda z)|\le \frac1{(1+\varepsilon)^k}\|f\|_{A}.
$$
Therefore,
$$
\sum_{k=0}^\infty \left| \varphi_s\left(\frac{ d^kf(\cdot)}{k!}(z)\right)\right|\le \frac{1+\varepsilon}{\varepsilon}\|f\|_{A},
$$
which shows that $\varphi_s^z$ is a well  defined and continuous homomorphism, that is, an element of $M_b(B_{\ell_p})$. Hence, $\dist_{M_{b}(B_{\ell_p})}(\varphi_s)\ge\big(1-s^p\big)^{1/p}$.

Now,  $1-s$ goes to 0 faster than $\big(1-s^p\big)^{1/p}$ as $s\to 1^-$, so the result follows.
$\square$ \end{proof}

\medskip

In the proof of the previous Proposition we have shown that $S({\varphi_s})$ contains the set
$$\{\varphi_s^z: z\in (1-s^p)^{1/p}B_{\ell_p}\},$$which can be thought of as a ball around $\varphi_s$.
For $p\in\mathbb N$, Proposition~\ref{hojas de Mb0} below shows that these two sets actually coincide.
This means that for any $\varphi\in M_b(B_{\ell_p})$ defined as in the Proposition, the sheet of $\varphi$ is (via $\pi$) a copy of a centered ball of radius strictly smaller than one.
Therefore, $M_b(B_{\ell_p})$ cannot be seen as a union of disjoint copies (via $\pi$) of $B_{\ell_p}$, as one might have thought from the case $U=E$, where $M_b(E)$ is a disjoint union of analytic copies of $E''$.

\begin{definition}
For $0<r<1$, we will say that $\varphi$ is an \emph{$r$-block homomorphism} on $H_b(B_{\ell_p})$ if $\varphi$ is an adherent point in $M_b(B_{\ell_p})$ of evaluations at the elements of a block basic sequence  $(x_n)_{n\in\mathbb N}$ with $\|x_n\|=r$.
\end{definition}
Since block bases are weakly null, any $r$-block homomorphism belongs to $\pi^{-1}(0)$. Also, adherent points of a block basic sequence  $(y_n)_{n\in\mathbb N}$ in the unit ball of $\ell_p$ with $\|y_n\|\to r$ are $r$-block homomorphisms, since they are also adherent points of $(ry_n/\|y_n\|)_{n\in\mathbb N}$.
Now we describe the sheet of an $r$-block homomorphism.

\begin{proposition}\label{hojas de Mb0}
Let $p$ be a natural number greater than 1 and, for $0<r<1$,  let $\varphi$ be an $r$-block homomorphism on $H_b(B_{\ell_p})$. Then $$S(\varphi)=\{\varphi^z:  \|z\|^p+r^p< 1 \}.$$
\end{proposition}
\begin{proof}We take $(x_n)_n$ a block basic sequence with $\|x_n\|=r$ such that $\varphi$ is an adherent point of $(\delta_{x_n})_n$. Proceeding as in the proof of Proposition~\ref{no reverse ineq} (with $x_n$ in the role of $s e_n$), we can show that $\varphi^z$ is a well defined homomorphism for $\|z\|^p<1-r^p$. This gives one inclusion.

For the reverse inclusion, $(x_n)_n$ is a block basic sequence, so there exists a sequence $(\alpha_k)_k\subset\mathbb C$ and finite subsets $J_n$ of $\mathbb N$ with $\max J_n< \min J_{n+1}$ such that $$x_n=\sum_{k\in J_n}\alpha_ke_k \text{ for all }n\ge 1.$$
Suppose $\|z\|^p>1-r^p$. Since $\|x_n\|=r$ for every $n\in\mathbb N$, there are some $\delta>0$ and $M\in\mathbb N$ such that for all $n>M$,
$$
\Big|\sum_{k=1}^M|z_k|^p + \sum_{k\in J_n}\big|z_k+\alpha_k\big|^p -\sum_{
\overset{k>M}{^{k\notin J_n}}
}|z_k|^p\Big|>1+\delta.
$$
Let us define $f_{N}(x)=\Big(\sum_{k=1}^\infty\theta_kx_k^p\Big)^N$, where $\theta_k$ is a modulus one complex number with $\theta_kz_k^p=|z_k|^p$ for $1\le k\le M$, $\theta_k(z_k+\alpha_k)^p=|z_k+\alpha_k|^p$ for  $k\in J_n$ with $n>M$, and $\theta_k=1$ otherwise. The sequence  $\{f_{N}:N\in\mathbb N\}$ is bounded in $H_b(B_{\ell_p})$, since it consists of norm one homogeneous polynomials. Also, for $n>M$ we have
$$
|f_{N}(z+x_n)|\ge\Big|\sum_{k=1}^M|z_k|^p +\sum_{k\in J_n}\big|\alpha_k+z_k\big|^p -\sum_{
\overset{k>M}{^{k\notin J_n}}
}|z_k|^p\Big|^N>(1+\delta)^N.
$$
Since $f_N$ is a polynomial, $\varphi^z(f_{N})=\varphi( f_{N}(z+ \cdot))$, which is a limit point of $(f_N(z+x_n))_n$.
Therefore, $|\varphi^z(f_{N})|\ge(1+\delta)^N$, which implies that $\varphi^z$ cannot be a continuous homomorphism on $H_b(B_{\ell_p})$.
We have shown that
$S(\varphi)\subset \{\varphi^z:  \|z\|^p+r^p\le 1 \}.$ Since $S(\varphi)$ is open and $\pi$ is a local homeomorphism, we have the desired inclusion.
$\square$ \end{proof}

Now we  show that if we restrict ourselves to a distinguished part of the spectrum of $B_{\ell_p}$, with $p$ a natural number greater than 1, then the extension of every function in $H_b(B_{\ell_p})$ is of bounded type.
Let us define the subdomain $M^0_b(B_{\ell_p})$ of  $M_b(B_{\ell_p})$ as the union of the sheets of all $r$-block homomorphisms.
All adherent points of the sequence $(\delta_{r e_n})_n$ with $0<r<1$ are $r$-block homomorphisms, so the number of connected components of $
M^0_b(B_{\ell_p})$ has at least the cardinality of $\beta\mathbb N$. Moreover, it is not clear that there are morphisms in $M_b(B_{\ell_p})$ that are not in $M^0_b(B_{\ell_p})$ (though to assert such a thing one should be able to prove of a really strong Corona theorem for $H_b(B_{\ell_p})$). One might argue that morphisms in $M_b(B_{\ell_p})$ can be built from sequences that are not blocks, or even with nets, but we cannot ensure that those do not have an alternative representation as $r$-block homomorphisms.  Anyway, $M^0_b(B_{\ell_p})$ is a relatively large part of $M_b(B_{\ell_p})$, where ``relatively'' should be understood as ``up to our knowledge".

Our distinguished spectrum $M_b^0(B_{\ell_p})$ is an open subset of $M_b(B_{\ell_p})$, since it is the union of some connected components of $M_b(B_{\ell_p})$. Thus $M_b^0(B_{\ell_p})$ is a Riemann domain over $\ell_p$ and every function $f\in H_b(B_{\ell_p})$ extends to a holomorphic function $\tilde f$ on $M_b^0(B_{\ell_p})$. Also, we have $\dist_{M_b^0(B_{\ell_p})}(\phi)=\dist_{M_b(B_{\ell_p})}(\phi)$ for any $\phi\in M_b^0(B_{\ell_p})$. We now show that this extension is of bounded type.

\begin{proposition}\label{extension al distinguido}
 Let $p$ be a natural number greater than 1, then the canonical extension $\tilde f$ of any  $f\in H_b(B_{\ell_p})$ to $M_b^0(B_{\ell_p})$ is of bounded type.
\end{proposition}
\begin{proof}
Given $\varepsilon>0$ we consider the $M_b^0(B_{\ell_p})$-bounded set $$A=\{\phi\in M_b^0(B_{\ell_p}):\, \dist_{M_b^0(B_{\ell_p})}(\phi)>\varepsilon\}.$$ By Proposition \ref{hojas de Mb0}, $A$ intersects the sheet of an $r$-block homomorphism $\varphi$ if and only if $1-r^p>\varepsilon^p$. Moreover, for such $\varphi$ and $r$ we have $$A\cap S(\varphi)=\{\varphi^z:\; (\|z\|+\varepsilon)^p+r^p <1 \}.$$
Let $(x_n)_n$ be a block basic sequence with $\|x_n\|=r$ such that $\varphi$ is an adherent point of $(\delta_{x_n})_n$. Fix $z$ with $(\|z\|+\varepsilon)^p+r^p <1$. 
Since $(x_n)_n$ is a block sequence of elements of norm $r$, there exists $n_0$ such that for $n\ge n_0$ we get
\begin{equation}\label{desigualdades normas}\|z+x_n\|^p<\|z\|^p+\|x_n\|^p+\varepsilon^p/2 <(\|z\|+\varepsilon)^p+\|x_n\|^p-\varepsilon^p/2<1-\varepsilon^p/2.
\end{equation}
If we set $R=(1-\varepsilon^p/2)^{1/p}$, equation (\ref{desigualdades normas}) shows that $\{x_n+z:\, n\ge n_0\}$ is contained in $R\, B_{\ell_p}$. From this we  conclude that  $\varphi^z\prec R\, B_{\ell_p}$. Now, $R$ depends only on $\varepsilon$, and any  $\phi \in A$ is of the form $\varphi^z$ as above. Therefore, the extension $\tilde f$ to $M_b^0(B_{\ell_p})$ of any $f\in H_b(B_{\ell_p})$ satisfies
$\|\tilde f\|_A=\sup_{\phi\in A}|\phi(f)|\le \|f\|_{R\, B_{\ell_p}}<\infty$. This completes the proof.
$\square$ \end{proof}

The homomorphisms shown in Proposition~\ref{no reverse ineq} actually belong to $M_b^0(B_{\ell_p}),$ since they are adherent points of evaluations in multiple of elements of the canonical basis.  As a consequence, we see that there is no reverse inequality in (\ref{desigualdad dist}) even if we restrict ourselves to $M_b^0(B_{\ell_p})$. On the other hand, Proposition~\ref{extension al distinguido} shows that extensions to $M_b^0(B_{\ell_p})$ are of bounded type. One may then think that, if there are functions in $H_b(B_{\ell_p})$ whose extensions to $M_b(B_{\ell_p})$  fail to be of bounded type, their existence will probably not be based on the absence of the reverse inequality in~(\ref{desigualdad dist}).

Finally, since by Proposition~\ref{hojas de Mb0} the connected components of $M_b^0(B_{\ell_p})$ are balls, we have the following corollary, which could also be deduced from the last proposition, proceeding as in the proof of \cite[Proposition 2.4]{DinVen04}:

\begin{corollary}
If $p$ is a natural number greater than 1, then $M_b^0(B_{\ell_p})$ is an $H_b$-domain of holomorphy.
\end{corollary}

\end{document}